\documentclass{amsart}

\usepackage{amsmath,amssymb,amsthm}
\usepackage{enumerate}
\usepackage[all]{xy}

\theoremstyle{plain}
\newtheorem{theorem}{Theorem}[section]
\newtheorem*{theorem*}{Theorem}
\newtheorem{lemma}[theorem]{Lemma}

\newtheorem{corollary}[theorem]{Corollary}

\theoremstyle{definition}

\theoremstyle{remark}
\newtheorem{remark}[theorem]{Remark}
\newtheorem{example}[theorem]{Example}

\DeclareMathOperator\Aut{Aut}
\DeclareMathOperator\supp{supp}
\DeclareMathOperator\disp{disp}
\DeclareMathOperator\diam{diam}
\newcommand\id{\mathrm{id}}

\newcommand\h[1]{\nobreakdash-\hspace{0pt}} % no break at hyphen,
\newcommand\C{\ensuremath{\mathcal{C}}}
\newcommand\M{\ensuremath{\mathcal{M}}}
\renewcommand\H{\ensuremath{\mathcal{H}}}
\let\setminus\smallsetminus

\newcommand\Cbar{\overline{C}}
\def\<#1>{\langle #1\rangle}

\begin{document}

\title[Automorphisms of buildings]{Automorphisms of Non-Spherical Buildings\\ Have
  Unbounded Displacement}

\author[Abramenko]{Peter Abramenko}
\address{Department of Mathematics\\
University of Virginia\\
Charlottesville, VA 22904}
\email{pa8e@virginia.edu}

\author[Brown]{Kenneth S. Brown}
\address{Department of Mathematics\\ 
Cornell University\\ 
Ithaca, NY 14853}
\email{kbrown@cornell.edu}

\date{October 5, 2007}

\begin{abstract}
If $\phi$ is a nontrivial automorphism of a thick building~$\Delta$ of
purely infinite type, we prove that there is no bound on the distance
that $\phi$ moves a chamber.  This has the following group-theoretic
consequence:  If $G$ is a group of automorphisms of~$\Delta$ with
bounded quotient, then the center of~$G$ is trivial.
\end{abstract}

\maketitle

\section*{Introduction}
\label{sec:introduction}

A well-known folklore result says that a nontrivial
automorphism~$\phi$ of a thick Euclidean building~$X$ has unbounded
displacement.  Here we are thinking of $X$ as a metric space, and the
assertion is that there is no bound on the distance that $\phi$ moves
a point.  [For the proof, consider the action of~$\phi$ on the
boundary~$X_\infty$ at infinity.  If $\phi$ had bounded displacement,
then $\phi$ would act as the identity on~$X_\infty$, and one would
easily conclude that $\phi=\id$.]  In this note we generalize this
result to buildings that are not necessarily Euclidean.  We work with
buildings~$\Delta$ as combinatorial objects, whose set $\C$ of
chambers has a discrete metric (``gallery distance'').  We say that
$\Delta$ is of \emph{purely infinite type} if every irreducible factor
of its Weyl group is infinite.

\begin{theorem*}
Let $\phi$ be a nontrivial automorphism of a thick building~$\Delta$
of purely infinite type.  Then $\phi$, viewed as an isometry of the
set \C\ of chambers, has unbounded displacement.
\end{theorem*}

The crux of the proof is a result about Coxeter groups (Lemma
\ref{lem:3}) that may be of independent interest.  We prove the lemma
in Section~\ref{sec:lemma-about-coxeter}, after a review of the Tits
cone in Section~\ref{sec:preliminaries}.  We then prove the theorem in
Section~\ref{sec:proof-theorem}, and we obtain the following (almost
immediate) corollary:  If $G$ is a subgroup of~$\Aut(\Delta)$ such
that there is a bounded set of representatives for the $G$\h-orbits
in~\C, then the center of~$G$ is trivial.

We conclude the paper with a brief discussion of displacement in the
spherical case.  We are grateful to Hendrik Van Maldeghem for
providing us with some counterexamples in this connection (see
Example~\ref{exam:1} and Remark~\ref{rem:5}).

\section{Preliminaries on the Tits cone}
\label{sec:preliminaries}

In this section we review some facts about the Tits cone associated to
a Coxeter group
\cite{abramenko08:_approac_to_build,bourbaki81:_group_lie,humphreys90:_reflec_coxet,tits61:_group_coxet,vinberg71:_discr}.
We will use \cite{abramenko08:_approac_to_build} as our basic
reference, but much of what we say can also be found in one or more of
the other cited references.

Let $(W,S)$ be a Coxeter system with $S$ finite.  Then $W$ admits a
canonical representation, which turns out to be faithful (see
Lemma~\ref{lem:7} below), as a linear reflection group acting on a real
vector space $V$ with a basis $\{e_s \mid s\in S\}$.  There is an
induced action of~$W$ on the dual space~$V^*$.  We denote by~$C_0$ the
simplicial cone in~$V^*$ defined by
\[
C_0 := \{x\in V^* \mid \<x,e_s> >0 \text{ for all } s\in S\};
\]
here $\<-,->$ denotes the canonical evaluation pairing between $V^*$
and~$V$.  We call~$C_0$ the \emph{fundamental chamber}.  For each
subset $J\subseteq S$, we set
\[
A_J := \{x \in V^* \mid \<x,e_s> = 0 \text{ for } s \in J \text{
  and } \<x,e_s> > 0 \text{ for } s \in S\setminus J\}.
\]
The sets $A_J$ are the (relatively open) \emph{faces} of~$C$ in the
standard terminology of polyhedral geometry.  They form a partition of
the closure~$\Cbar_0$ of $C_0$ in~$V^*$.

For each $s\in S$, we denote by~$H_s$ the hyperplane in~$V^*$ defined
by the linear equation $\<-,e_s> =0$.  If follows from the explicit
definition of the canonical representation of~$W$ (which we have not
given) that $H_s$ is the fixed hyperplane of $s$ acting on~$V^*$.  The
complement of $H_s$ in~$V^*$ is the union of two open
halfspaces~$U_\pm(s)$ that are interchanged by~$s$.  Here
\[
U_+(s) := \{x \in V^* \mid \<x,e_s> > 0\},
\]
and
\[
U_-(s) := \{x \in V^* \mid \<x,e_s> < 0\}.
\]
The hyperplanes $H_s$ are called the \emph{walls} of~$C_0$.  We denote
by $\H_0$ the set of walls of~$C_0$.

The \emph{support} of the face
$A= A_J$, denoted~$\supp A$, is defined to be the intersection of the
walls of~$C_0$ containing~$A$, i.e., $\supp A = \bigcap_{s\in J} H_s$.
Note that $A$ is open in~$\supp A$ and that $\supp A$ is the linear
span of~$A$.

Although our definitions above made use of the basis $\{e_s \mid s \in
S\}$ of~$V$, there are also intrinsic geometric characterizations of
walls and faces.  Namely, the walls of~$C_0$ are the hyperplanes $H$
in~$V^*$ such that $H$ does not meet~$C_0$ and $H\cap\Cbar_0$ has nonempty
interior in~$H$.  And the faces of~$C_0$ correspond to subsets
$\H_1\subseteq \H_0$.  Given such a subset, let $L :=
\bigcap_{H\in\H_1} H$; the corresponding face~$A$ is then the relative
interior (in~$L$) of the intersection $L\cap\Cbar_0$.

We now make everything $W$-equivariant.  We call a subset $C$ of~$V^*$
a \emph{chamber} if it is of the form $C = wC_0$ for some $w \in W$,
and we call a subset $A$ of~$V^*$ a \emph{cell} if it is of the form
$A = wA_J$ for some $w \in W$ and $J \subseteq S$.  Each chamber~$C$
is a simplicial cone and hence has well-defined walls and faces, which
can be characterized intrinsically as above.  If $C=wC_0$ with $w\in
W$, the walls of~$C$ are the transforms $wH_s$ ($s \in S$), and the
faces of~$C$ are the cells $wA_J$ ($J\subseteq S$).  Finally, we call
a hyperplane $H$ in~$V^*$ a \emph{wall} if it is a wall of some
chamber, and we denote by~\H\ the set of all walls; thus
\[
\H = \{wH_s \mid w \in W,\, s \in S\}.
\]
The set of all faces of all chambers is equal to the set of all cells.
The union of these cells is called the \emph{Tits cone} and will be
denoted by~$X$ in the following.  Equivalently,
\[
X = \bigcup_{w\in W} w\Cbar_0.
\]

We now record, for ease of reference, some standard facts about the
Tits cone.  The first fact is Lemma~2.58 in
\cite[Section~2.5]{abramenko08:_approac_to_build}.  See also the proof
of Theorem~1 in \cite[Section~V.4.4]{bourbaki81:_group_lie}.

\begin{lemma}
\label{lem:6}
For any $w \in W$ and $s \in S$, we have
\[
wC_0 \subseteq U_+(s) \iff l(sw) > l(w)
\]
and
\[
wC_0 \subseteq U_-(s) \iff l(sw) < l(w).
\]
Here $l(-)$ is the length function on~$W$ with respect to~$S$.
\qed
\end{lemma}

This immediately implies:

\begin{lemma}
\label{lem:7}
$W$ acts simply transitively on the set of chambers.
\qed
\end{lemma}

The next result allows one to talk about separation of cells by
walls.  It is part of Theorem~2.80 in
\cite[Section~2.6]{abramenko08:_approac_to_build}, and it can also be
deduced from Proposition~5 in \cite[Section~V.4.6]{bourbaki81:_group_lie}.

\begin{lemma}
\label{lem:8}
If $H$ is a wall and $A$ is a cell, then either $A$ is contained
in~$H$ or $A$ is contained in one of the two open halfspaces
determined by~$H$.  \qed
\end{lemma}

We turn now to reflections.  The following lemma is an easy
consequence of the stabilizer calculation in
\cite[Theorem~2.80]{abramenko08:_approac_to_build} or
\cite[Section~V.4.6]{bourbaki81:_group_lie}.

\begin{lemma}
\label{lem:9}
For each wall $H\in\H$, there is a unique nontrivial element $s_H\in
W$ that fixes $H$ pointwise.
\qed
\end{lemma}
We call $s_H$ the \emph{reflection} with respect to~$H$.  In view of a
fact stated above, we have $s_{H_s} = s$ for all $s\in S$.  Thus $S$
is the set of reflections with respect to the walls in~$\H_0$.  It
follows immediately from Lemma~\ref{lem:9} that
\begin{equation}
\label{eq:1}
s_{wH} = ws_Hw^{-1}
\end{equation}
for all $H \in \H$ and $w \in W$.  Hence $wSw^{-1}$ is the set of
reflections with respect to the walls of~$wC_0$.

\begin{corollary}
\label{cor:2}
For $s \in S$ and $w \in W$, $\,H_s$ is a wall of~$wC_0$ if and only
if $w^{-1}sw$ is in~$S$.
\end{corollary}

\begin{proof}
$H_s$  is a wall of~$wC_0$ if and only if $s$ is the reflection with
respect to a wall of~$wC_0$.  In view of the observations above, this
is equivalent to saying $s\in wSw^{-1}$, i.e., $w^{-1} s w \in S$.
\end{proof}

Finally, we record some special features of the infinite case.

\begin{lemma}
\label{lem:1}
Assume that $(W,S)$ is irreducible and $W$ is infinite.
\begin{enumerate}
\item\label{item:1}
If two chambers $C,D$ have the same walls, then $C=D$.
\item\label{item:2}
The Tits cone $X$ does not contain any pair~$\pm x$ of opposite
nonzero vectors.
\end{enumerate}
\end{lemma}

\begin{proof}
(\ref{item:1}) We may assume that $C=C_0$ and~$D =wC_0$ for some $w\in
W$.  Then Corollary~\ref{cor:2} implies that $C$ and~$D$ have the same
walls if and only if $w$ normalizes~$S$.  So the content
of~(\ref{item:1}) is that the normalizer of $S$ in~$W$ is trivial.  This
is a well known fact.  See \cite[Section~V.4,
Exercise~3]{bourbaki81:_group_lie},
\cite[Proposition~4.1]{deodhar82:_coxet}, or
\cite[Section~2.5.6]{abramenko08:_approac_to_build}.
Alternatively, there is a direct geometric proof of~(\ref{item:1})
outlined in the solution to
\cite[Exercise~3.118]{abramenko08:_approac_to_build}.

(\ref{item:2}) This is a result of Vinberg~\cite[p.~1112,
Lemma~15]{vinberg71:_discr}.  See also
\cite[Section~2.6.3]{abramenko08:_approac_to_build} and
\cite[Theorem~2.1.6]{krammer94:_coxet} for alternate proofs.
\end{proof}

\section{A lemma about Coxeter groups}
\label{sec:lemma-about-coxeter}

We begin with a geometric version of our lemma, and then we translate
it into algebraic language.

\begin{lemma}
\label{lem:2}
Let $(W,S)$ be an infinite irreducible Coxeter system with $S$ finite.
If $C$ and~$D$ are distinct chambers in the Tits cone, then $C$ has a
wall~$H$ with the following two properties:
\begin{enumerate}[\rm(a)]
\item
$H$ is not a wall of~$D$.
\item
$H$ does not separate $C$ from~$D$.
\end{enumerate}
\end{lemma}

\begin{proof}
For convenience (and without loss of generality), we assume that $C$
is the fundamental chamber~$C_0$.  Define $J \subseteq S$ by
\[
J := \{s \in S \mid H_s \text{ is a wall of } D\},
\]
and set $L := \bigcap_{s\in J} H_s$.  Thus $L$ is the support of
the face $A = A_J$ of~$C$.  By Lemma~\ref{lem:1}(\ref{item:1}), $J\neq
S$, hence $L \neq \{0\}$.  Since $L$ is an intersection of walls
of~$D$, it is also the support of a face $B$ of~$D$.  Note that $A$
and~$B$ are contained in precisely the same walls, since they have the
same span~$L$.  In particular, $B$ is not contained in any of the
walls~$H_s$ with $s\in S\setminus J$, so, by Lemma~\ref{lem:8}, $B$ is
contained in either $U_+(s)$ or~$U_-(s)$ for each such~$s$.

Suppose that $B\subseteq U_-(s)$ for each $s \in S\setminus J$.  Then,
in view of the definition of~$A = A_J$ by linear equalities and
inequalities, $B \subseteq -A$.  But $B$ contains a nonzero vector~$x$
(since $B$ spans~$L$), so we have contradicted
Lemma~\ref{lem:1}(\ref{item:2}).  Thus there must exist $s\in S\setminus
J$ with $B \subseteq U_+(s)$.  This implies that $D\subseteq U_+(s)$,
and the wall $H=H_s$ then has the desired properties (a) and~(b).
\end{proof}

We now prove the algebraic version of the lemma, for which we relax
the hypotheses slightly.  We do not even have to assume that $S$ is
finite.  Recall that $(W,S)$ is said to be \emph{purely infinite} if
each of its irreducible factors is infinite.

\begin{lemma}
\label{lem:3}
Let $(W,S)$ be a purely infinite Coxeter system.  If $w \neq 1$
in~$W$, then there exists $s\in S$ such that:
\begin{enumerate}[\rm(a)]
\item
$w^{-1} s w \notin S$.
\item
$l(sw) > l(w)$.
\end{enumerate}
\end{lemma}

\begin{proof}
Let $(W_i,S_i)$ be the irreducible factors of~$(W,S)$, which are all
infinite.  Suppose the lemma is true for each factor~$(W_i,S_i)$, and
consider any $w\neq 1$ in~$W$.  Then $w$ has components $w_i\in W_i$,
at least one of which (say~$w_1$) is nontrivial.  So we can find $s
\in S_1$ with $w_1^{-1}sw_1 \notin S_1$ and $l(sw_1) > l(w_1)$.  One
easily deduces (a) and~(b).  We are now reduced to the case where
$(W,S)$ is irreducible.

If $S$ is finite, we apply Lemma~\ref{lem:2} with $C$ equal to the
fundamental chamber~$C_0$ and $D = w C_0$.  Then $H = H_s$ for some
$s\in S$.  Property~(a) of that lemma translates to (a) of the present
lemma by Corollary~\ref{cor:2}, and property~(b) of that lemma
translates to (b) of the present lemma by Lemma~\ref{lem:6}.

If $S$ is infinite, we use a completely different method.  The result
in this case follows from Lemma~\ref{lem:4} below.
\end{proof}

Recall that for any Coxeter system $(W,S)$ and any $w\in W$, there is
a (finite) subset $S(w) \subseteq S$ such that every reduced
decomposition of~$w$ involves precisely the generators in~$S(w)$.
This follows, for example, from Tits's solution to the word
problem~\cite{tits69:_word_prob}.  (See also
\cite[Theorem~2.35]{abramenko08:_approac_to_build}).

\begin{lemma}
\label{lem:4}
Let $(W,S)$ be an irreducible Coxeter system, and let $w \in W$ be
nontrivial.  If $S(w) \neq S$, then there exists $s\in S$ satisfying
conditions (a) and~(b) of Lemma~\ref{lem:3}.
\end{lemma}

\begin{proof}
By irreducibility, there exists an element $s\in S \setminus S(w)$
that does not commute with all elements of~$S(w)$.  Condition~(b) then
follows from the fact that $s\notin S(w)$ and standard properties of
Coxeter groups; see \cite[Lemma~2.15]{abramenko08:_approac_to_build}.
To prove~(a), suppose $sw=wt$ with $t\in S$.  We have $s\notin S(w)$ but
$s \in S(sw)$ (since $l(sw) > l(w)$), so necessarily $t = s$.
Using induction on~$l(w)$, one now deduces from Tits's solution to the
word problem that $s$ commutes with every element of~$S(w)$ (see
\cite[Lemma~2.39]{abramenko08:_approac_to_build}), contradicting the
choice of~$s$.
\end{proof}

\section{Proof of the theorem}
\label{sec:proof-theorem}

In this section we assume familiarity with basic concepts from the
theory of buildings
\cite{abramenko08:_approac_to_build,ronan89:_lectur,scharlau95:_build,tits74:_build_bn,weiss03}.

Let $\Delta$ be a building with Weyl group~$(W,S)$, let \C\ be the set
of chambers of~$\Delta$, and let $\delta\colon \C\times\C\to W$ be the
Weyl distance function.  (See
\cite[Section 4.8 or~5.1]{abramenko08:_approac_to_build} for the
definition and standard properties of~$\delta$.)  Recall that \C\ has
a natural \emph{gallery metric}~$d(-,-)$ and that
\[
d(C,D) = l\bigl(\delta(C,D)\bigr)
\]
for $C,D\in\C$.  Let $\phi \colon \Delta\to\Delta$ be an automorphism
of~$\Delta$ that is not necessarily type-preserving.  Recall that
$\phi$ induces an automorphism $\sigma$ of~$(W,S)$.  From the
simplicial point of view, we can think of $\sigma$ (restricted to~$S$)
as describing the effect of~$\phi$ on types of vertices.  From the
point of view of Weyl distance, $\sigma$ is characterized by the
equation
\[
\delta(\phi(C),\phi(D)) = \sigma(\delta(C,D))
\]
for $C,D\in\C$.

Our main theorem will be obtained from the following technical lemma:

\begin{lemma}
\label{lem:5}
Assume that $\Delta$ is thick.  Fix a chamber $C \in \C$, and set $w
:= \delta(C,\phi(C))$.  Suppose there exists $s \in S$ such that
$l(sw) > l(w)$ and $w^{-1}sw \neq t := \sigma(s)$.  Then there is a
chamber $D$ $s$\h-adjacent to~$C$ such that $d(D,\phi(D)) > d(C,
\phi(C))$.
\end{lemma}

\begin{proof}
We will choose $D$ $s$\h-adjacent to~$C$ so that $u :=
\delta(C,\phi(D))$ satisfies $l(u) \geq l(w)$ and $l(su) > l(u)$.  We
will then have $\delta(D,\phi(D)) = su$, as illustrated in the
following schematic diagram:
\[
\def\labelstyle{\textstyle}
\xymatrix@+1pc{C \ar[r]^w \ar@{-}[d]_s \ar[dr]_u & \phi(C) \ar@{-}[d]^t\\
  D \ar[r]_{su} & \phi(D)}
\]
Hence
\[
d(D,\phi(D)) = l(su) > l(u) \geq l(w) = d(C,\phi(C)).
\]

Case 1.  $l(wt) < l(w)$.  Then there is a unique chamber $E_0$
$t$\h-adjacent to~$\phi(C)$ such that $\delta(C,E_0) = wt$.  For all
other $E$ that are $t$\h-adjacent to~$\phi(C)$, we have $\delta(C,E) =
w$.  So we need only choose $D$ so that $\phi(D) \neq E_0$; this is
possible by thickness.  Then $u = w$.

Case 2.  $l(wt) > l(w)$.  Then $l(swt) > l(wt)$ because the conditions
$l(swt) < l(wt)$, $\,l(wt) > l(w)$, and $l(sw) > l(w)$ would imply
(e.g., by the deletion condition for Coxeter groups) $swt = w$, and
the latter is excluded by assumption.  Hence in this case we can
choose $D$ $s$\h-adjacent to~$C$ arbitrarily, and we then have $u =
wt$.
\end{proof}

Suppose now that $(W,S)$ is purely infinite and $\phi$ is nontrivial.
Then we can start with any chamber~$C$ such that $\phi(C)\neq C$, and
Lemma~\ref{lem:3} shows that the hypothesis of Lemma~\ref{lem:5} is
satisfied.  We therefore obtain a chamber~$D$ such that $d(D,\phi(D))
> d(C,\phi(C))$.  Our main theorem as stated in the introduction
follows at once.  We restate it here for ease of reference:

\begin{theorem}
\label{thr:1}
Let $\phi$ be a nontrivial automorphism of a thick building~$\Delta$
of purely infinite type.  Then $\phi$, viewed as an isometry of the
set \C\ of chambers, has unbounded displacement, i.e., the set
$\{d(C,\phi(C)) \mid C\in\C\}$ is unbounded.
\qed
\end{theorem}

\begin{remark}
\label{rem:1}
Note that, in view of the generality under which we proved Lemma
\ref{lem:3}, the building~$\Delta$ is allowed to have infinite rank.
\end{remark}

\begin{remark}
\label{rem:2}
In view of the existence of translations in Euclidean Coxeter
complexes, the thickness assumption in the theorem cannot be dropped.
\end{remark}

\begin{corollary}
\label{cor:1}
Let $\Delta$ and~\C\ be as in the theorem, and let $G$ be a group of
automorphisms of~$\Delta$.  If there is a bounded set of
representatives for the $G$\h-orbits in~\C, then $G$ has trivial center.
\end{corollary}

\begin{proof}
Let \M\ be a bounded set of representatives for the $G$\h-orbits
in~\C, and let $z\in G$ be central.  Then there is an upper bound~$M$
on the distances $d(C,zC)$ for $C\in\M$; we can take $M$ to be the
diameter of the bounded set $\M \cup z\M$, for instance.  Now every
chamber $D\in\C$ has the form $D=gC$ for some $g\in G$ and $C\in\M$,
hence
\[
d(D,zD) = d(gC,zgC) = d(gC,gzC) = d(C,zC) \leq M.
\]
Thus $z$ has bounded displacement and therefore $z=1$ by the theorem.
\end{proof}

\begin{remark}
\label{rem:4}
Although Corollary~\ref{cor:1} is stated for faithful group actions,
we can also apply it to actions that are not necessarily faithful and
conclude (under the hypothesis of the corollary) that the center
of~$G$ acts trivially.
\end{remark}

\begin{remark}
\label{rem:6}
Note that the hypothesis of the corollary is satisfied if the action
of $G$ is chamber transitive.  In particular, it is satisfied if the
action is strongly transitive and hence corresponds to a BN-pair
in~$G$.  In this case, however, the result is trivial (and does not
require the building to be of purely infinite type).  Indeed, the
stabilizer of every chamber is a parabolic subgroup and hence is
self-normalizing, so it automatically contains the center of~$G$.  To
obtain other examples, consider a cocompact action of a group on a
locally-finite thick Euclidean building (e.g., a thick tree).  The
corollary then implies that the center of the group must act
trivially.
\end{remark}

\begin{remark}
\label{rem:3}
The conclusion of Theorem~\ref{thr:1} is obviously false for spherical
buildings, since the metric space~\C\ is bounded in this case.  But
one can ask instead whether or not
\begin{equation}
\label{eq:2}
\disp \phi = \diam \Delta,
\end{equation}
where $\diam \Delta$ denotes the diameter of the metric space~\C, and
$\disp\phi$ is the \emph{displacement} of~$\phi$; the latter is defined by
\[
\disp \phi := \sup \{d(C,\phi(C)) \mid C\in\C\}.
\]
Note that, in the spherical case, Equation~\eqref{eq:2} holds if and
only if there is a chamber~$C$ such that $\phi(C)$ and~$C$ are
opposite.  This turns out to be false in general.  The following
counterexample was pointed out to us by Hendrik Van Maldeghem.
\end{remark}

\begin{example}
\label{exam:1}
Let $k$ be a field and $n$ an integer $\geq 2$.  Let $\Delta$ be the
building associated to the vector space $V=k^{2n}$.  Thus the vertices
of~$\Delta$ are the subspaces $U$ of~$V$ such that $0<U<V$, and the
simplices are the chains of such subspaces.  A chamber is a chain
\[
U_1<U_2<\cdots <U_{2n-1}
\]
with $\dim U_i = i$ for all~$i$, and two such chambers $(U_i)$
and~$(U'_i)$ are opposite if and only if $U_i + U'_{2n-i} = V$ for
all~$i$.  Now choose a non-degenerate alternating bilinear form $B$
on~$V$, and let $\phi$ be the (type-reversing) involution of~$\Delta$
that sends each vertex $U$ to its orthogonal subspace~$U^\perp$ with
respect to~$B$.  For any chamber $(U_i)$ as above, its image
under~$\phi$ is the chamber~$(U'_i)$ with $U'_{2n-i} = U_i^\perp$
for all~$i$.  Since $U_1\leq U_1^\perp = U'_{2n-1}$, these two
chambers are not opposite.
\end{example}

Even though \eqref{eq:2} is false in general, one can still use
Lemma~\ref{lem:5} to obtain lower bounds on~$\disp\phi$.  Consider the
rank~2 case, for example.  Then $\Delta$ is a generalized $m$-gon for
some~$m$, its diameter is~$m$, and its Weyl group~$W$ is the dihedral
group of order~$2m$.  Lemma~\ref{lem:5} in this case yields the
following result:

\begin{corollary}
\label{cor:3}
Let $\phi$ be a nontrivial automorphism of a thick generalized
$m$-gon.  Then the following hold:
\begin{enumerate}[\rm(a)]
\item
$\disp\phi \geq m-1$.
\item
If $\phi$ is type preserving and $m$ is odd, or if $\phi$ is type
reversing and $m$ is even, then $\disp\phi=m$.
\end{enumerate}
\end{corollary}

\begin{proof}
(a) The hypothesis of Lemma~\ref{lem:5} is always satisfied as long as
$w \neq 1$ and $l(w) < m-1$ since then there exists $s \in S$ with
$l(sw) > l(w)$, and $sw$ has a unique reduced decomposition (which is
in particular not of the form $wt$ with $t \in S$).

(b) Suppose $l(w) = m-1$, and let $s \in S$ be the unique element such
that $sw = w_0$, where the latter is the longest element of~$W$.  Then
$s' \in S$ satisfies $ws' = sw$ if and only if $l(ws') > l(w)$.  Since
$w$ does not start with~$s$, this is equivalent to $s' \neq s$ if $m$
is odd and to $s' = s$ if $m$ is even.  So the hypothesis of
Lemma~\ref{lem:5} is satisfied for \emph{any} $w$ different from 1
and~$w_0$ if $m$ is odd and $\sigma=\id$ or if $m$ is even and $\sigma
\neq \id$.
\end{proof}

We conclude by mentioning another family of examples, again pointed
out to us by Van Maldeghem.

\begin{remark}
\label{rem:5}
For even $m = 2n$, type-preserving automorphisms $\phi$ of generalized
$m$-gons with $\disp \phi = m-1$ arise as follows.  Assume that there
exists a vertex $x$ in the generalized $m$-gon~$\Delta$ such that the
ball $B(x, n)$ is fixed pointwise by~$\phi$.  Here $B(x,n)$ is the set
of vertices with $d(x,y) \leq n$, where $d(-,-)$ now denotes the usual
graph metric, obtained by minimizing lengths of paths.  Recall that
there are two types of vertices in~$\Delta$ and that opposite vertices
always have the same type since $m$ is even.  Let $y$ be any vertex
that does has not the same type as~$x$.  Then $y$ is at distance at
most $n-1$ from some vertex in $B(x,n)$.  Since $\phi$ fixes $B(x,n)$
pointwise, $d(y, \phi(y)) \leq 2n-2$.  So $C$ and~$\phi(C)$ are not
opposite for any chamber~$C$ having $y$ as a vertex.  Since this is
true for any vertex $y$ that does not have the same type as~$x$,
$\disp \phi \neq m$ and hence, by Corollary~\ref{cor:3}(a), $\disp \phi
= m-1$ if $\phi \neq \id$.  Now it is a well-known fact (see for
instance \cite[Corollary 5.4.7]{maldeghem98:_gener}) that every
\emph{Moufang} $m$-gon possesses nontrivial type-preserving
automorphisms $\phi$ fixing some ball $B(x,n)$ pointwise.  (In the
language of incidence geometry, these automorphisms are called central
or axial collineations, depending on whether $x$ is a point or a line
in the corresponding rank~2 geometry.)  So for $m = 4$, $6$, or~$8$,
all Moufang $m$-gons admit type-preserving automorphisms $\phi$ with
$\disp \phi = m-1$.
\end{remark}

\providecommand{\bysame}{\leavevmode\hbox to3em{\hrulefill}\thinspace}


\begin{thebibliography}{10}

\bibitem{abramenko08:_approac_to_build}
P.~Abramenko and K.~S. Brown, \emph{Approaches to buildings}, Springer, New
  York, 2008, in preparation.

\bibitem{bourbaki81:_group_lie}
N.~Bourbaki, \emph{\'{E}l\'ements de math\'ematique. {F}asc. {XXXIV}. {G}roupes
  et alg\`ebres de {L}ie. {C}hapitre {IV}: {G}roupes de {C}oxeter et syst\`emes
  de {T}its. {C}hapitre {V}: {G}roupes engendr\'es par des r\'eflexions.
  {C}hapitre {VI}: {S}yst\`emes de racines}, Actualit\'es Scientifiques et
  Industrielles, No. 1337, Hermann, Paris, 1968.

\bibitem{deodhar82:_coxet}
V.~V. Deodhar, \emph{On the root system of a {C}oxeter group}, Comm. Algebra
  \textbf{10} (1982), no.~6, 611--630.

\bibitem{humphreys90:_reflec_coxet}
J.~E. Humphreys, \emph{Reflection groups and {C}oxeter groups}, Cambridge
  Studies in Advanced Mathematics, vol.~29, Cambridge University Press,
  Cambridge, 1990.

\bibitem{krammer94:_coxet}
D.~Krammer, \emph{The conjugacy problem for {C}oxeter groups}, Ph.D. thesis,
  Universiteit Utrecht, 1994.

\bibitem{ronan89:_lectur}
M.~Ronan, \emph{Lectures on buildings}, Perspectives in Mathematics, vol.~7,
  Academic Press Inc., Boston, MA, 1989.

\bibitem{scharlau95:_build}
R.~Scharlau, \emph{Buildings}, Handbook of incidence geometry, North-Holland,
  Amsterdam, 1995, pp.~477--645.

\bibitem{tits61:_group_coxet}
J.~Tits, \emph{Groupes et g\'{e}om\'{e}tries de {C}oxeter}, Unpublished
  manuscript, 1961.

\bibitem{tits69:_word_prob}
\bysame, \emph{Le probl\`eme des mots dans les groupes de {C}oxeter}, Symposia
  {M}athematica ({INDAM}, {R}ome, 1967/68), Vol. 1, Academic Press, London,
  1969, pp.~175--185.

\bibitem{tits74:_build_bn}
\bysame, \emph{Buildings of spherical type and finite {BN}-pairs}, Lecture
  Notes in Mathematics, vol. 386, Springer-Verlag, Berlin, 1974.

\bibitem{maldeghem98:_gener}

H.~{V}an Maldeghem, \emph{Generalized polygons}, Monographs in Mathematics,
  vol.~93, Birkh\"auser Verlag, Basel, 1998.

\bibitem{vinberg71:_discr}
{\`E}.~B. Vinberg, \emph{Discrete linear groups that are generated by
  reflections}, Izv. Akad. Nauk SSSR Ser. Mat. \textbf{35} (1971), 1072--1112.

\bibitem{weiss03}
R.~M. Weiss, \emph{The structure of spherical buildings}, Princeton University
  Press, Princeton, NJ, 2003.

\end{thebibliography}
\end{document}